\documentclass[12pt,a4paper]{amsart}
\usepackage{amsmath,graphicx,amscd}

\theoremstyle{plain}
\newtheorem{theorem}{Theorem}[section]
\newtheorem{lemma}[theorem]{Lemma}

\newtheorem{proposition}[theorem]{Proposition}

\theoremstyle{remark}

\theoremstyle{definition}

\def \C {\mathbf{C}}

\def \S {\mathbf{S}}

\def\RR{\mathcal{R}}

\numberwithin{equation}{section}
\def\SU{{\rm SU}(2,1)}
\def\SLC{{\rm SL}(2,\C)}
\def\PU{{\rm PU}(2,1)}
\def\bhs{\Sigma(p,q,r)}

\def\int{\text{\bf in}}

\def\SUn{{\rm SU}(2)}

\def\fg{\pi_1}

\def\diag{{\rm diag}}

\begin{document}
\title{On the $\SU$ representation space of the Brieskorn homology spheres}
\author {Vu The Khoi}
\address{Institute of Mathematics, Vietnam Academy of Sciences and Technology, 18 Hoang Quoc Viet Road, 10307, Hanoi, Vietnam}
\email{vtkhoi@math.ac.vn}
\thanks{The author was supported by a COE Postdoctoral Fellowship of University of Tokyo, the JSPS's Kakenhi Grant and the National Basic Research Program of Vietnam}

\begin{abstract}
{In this paper, we give a parameterization of the $\SU$ representation space of the Brieskorn homology spheres using the trace coordinates. As applications, 
we give an example which shows that the orbifold Toledo invariant in \cite{krebs}  does not distinguish the connected components of the $\PU$ representation space.}
\end{abstract}
\maketitle
\section{Introduction} 
\vskip0.1cm
Let $M$ be a manifold with the fundamental group $\fg(M)$ and  $G$ be a Lie group. The \textit{representation space} of $M,$ denoted by $\RR_G(M),$ is the space of representations from $\fg(M)$ into the Lie group $G,$ modulo conjugation: 
$$\RR_G(M):= \hom(\fg(M),G)/G.$$ 
We denote by $\RR^*_G(M)$ the subset of the representation space which consists of irreducible representations.
The representation space of  $3$-manifolds has been studied extensively in the case where $G=\SUn, {\rm SU}(3)$ or  $\SLC$ in connection with the Casson invariants and
 hyperbolic geometry (see \cite{boden,bhc,cs,fs,furuta,kk}).

Let us recall that $\SU$ is the special unitary group corresponding to the indefinite inner product 
$\langle Z,W\rangle_{2,1}=Z_1\overline{W}_1 + Z_2\overline{W}_2 - Z_3\overline{W}_3$
on $\C^3.$ The group $\PU$ is the quotient of $\SU$ by its center. 

In this paper we study the $\SU$ representation space $\RR_{\SU}(M).$
 The motivation for this study comes from complex hyperbolic geometry where $\RR_{\PU}(M)$ serves as the local model for the deformation space of spherical
CR structures on $M.$ For convenience, we will work with the group $\SU$  and then deduce results for the $\PU$ case.  

   Let $p,q,r$ be pairwise coprime positive integers, the \textit{Brieskorn homology sphere} $\bhs$ is defined to be the link of  singularity in $\C^3,$ that is :
$$ \bhs:=\{(x,y,z)|\ x^p+y^q+z^r=0 \}\ \cap \S^5_{\epsilon}.$$

It is well known that the fundamental group of $\bhs$ may be given as 
$$\fg(\bhs)= \langle x,y,z,h| \ h\ \text{central},\ x^ph^a=y^qh^b=z^rh^c=xyz=1
\rangle,$$
where $a,b,c$ are integers satisfying
$$
\frac{a}{p}+ \frac{b}{q}+\frac{c}{r}=\frac{1}{pqr}.
$$

In this paper, for simplicity, we will denote by $t_A$ the trace of a matrix $A$ and $[A,B]$ the commutator
$ABA^{-1}B^{-1}.$  
The notations $\Re$ and $\Im$ 
stand for the real and imaginary part of a complex numbers
respectively.
 Our main theorem shows that $\RR^*_{\SU}(M)$ can be parameterized by certain trace coordinates. 

{\bf Theorem 3.1} {\it Two irreducible representations $\rho, \rho':\fg(\bhs)\longrightarrow \SU$ are conjugate if and only if
the image under  $\rho$ and $\rho'$ of 
each $x,y,h$ are conjugate and  
satisfy the relations
$$t_{\rho(xy)}=t_{\rho'(xy)},
\quad
t_{\rho(x^{-1}y)}=t_{\rho'(x^{-1}y)}, 
\quad
\Im(t_{\rho([x,y])})=\Im(t_{\rho'([x,y])}).$$  }

The rest of this paper is organized as follows. In section 2 we study the trace identities for the free group of rank two. Using algebraic results about the invariant
ring of matrices, we are able to deduce the coordinates and relations for the $\SU$ representation space of the free group of rank two.
Section 3 is devoted to the proof of the main result. In this section, we also show how to find the constraint for the parameters of the representation spaces
in practice. Finally, in section 4, we apply our results to give explicit descriptions of the representation spaces of the Brieskorn homology
spheres $\Sigma(2,3,11)$ and $\Sigma(2,3,13).$        

\section{Trace calculus for free group of rank two} 
\vskip0.1cm
 We first recall some known results about  matrices in $\SU$. 
The reader should consult \cite{chen,goldman} for details. 
Let $V_{\_}$ and $V_0$ be
the two subsets of $\C^3$ defined by
$V_{\_}:=\{Z=(z_1, z_2, z_3)\in \C^3|\ 
\langle Z,Z
\rangle_{2,1}\ <0 \}$ and $V_0:=\{Z=(z_1, z_2, z_3)\in \C^3|\ 
\langle Z,Z
\rangle_{2,1}\ =0 \}$.  We denote by 
$P:\C^3\,\setminus\{0\}\rightarrow \C P^2$
 the canonical projection onto the complex projective space. Then $P(V_{\_})$ equipped with the Bergman metric is  the model of the complex hyperbolic space $H^2_{\C}$. The boundary $\partial H^2_{\C}$ in $\C P^2$ is $P(V_0\setminus \{0\}).$

The
elements of $\SU$ can be classified according to their action on the complex hyperbolic space $H^2_{\C}$ \cite{chen}.
Namely, a matrix is called {\it elliptic} if it has a fixed point in $H^2_{\C}.$ We call it {\it parabolic} if it has a unique fixed point in $\overline {H^2_{\C}}$ which
lies on $\partial H^2_{\C}.$ And finally, an element is called {\it loxodromic} if it 
has exactly two fixed points in $\overline {H^2_{\C}}$ which lie on $\partial H^2_{\C}.$ 

A classification of conjugacy classes of elements of $\SU$ can be found in \cite{chen}. In particular
it says that two elliptic elements are conjugate if and only if they have the same positive and negative class of eigenvalues
(counted with multiplicity).
An explanation of terminology should be added here: we say that an eigenvalue $\lambda$ of an elliptic element is of {\it positive type} 
(respectively {\it negative type}) if it has an $\lambda$-eigenvector $v$ such that $
\langle v,v
\rangle_{2,1}$ is positive (respectively negative). It has been shown that
every eigenvalue of an elliptic element has either positive or negative type.         

The next proposition gives several trace identities for a pair of matrices in $\SU.$ These identities will be crucial in getting a coordinate system on the 
representation space. 
\begin{proposition}\label{trid} 
Let $A$ and $B$ be a pair of matrices in $\SU $.
Then the following equations hold:

{\rm i)} $t_{A^{-1}} = \overline{t_A}$.

{\rm ii)} $t_{A^2} = t_A^2-2\overline{t_A}$.

{\rm iii)} $t_{A^3} = t_A^3-3|t_A|^2 + 3$.

{\rm iv)} $t_{A^2B} = t_At_{AB} - \overline{t_A}t_B + t_{A^{-1}B}$.

{\rm v)} $t_{A^2B^2}= t_At_Bt_{AB}-t_A^2\overline{t_B}+t_A\overline{t_{A^{-1}B}}-\overline{t_A}t_B^2+\overline{t_At_B}+t_Bt_{A^{-1}B}$.

{\rm vi)} $t_{ABAB^{-1}}= t_{AB}\overline{t_{A^{-1}B}}+\overline{t_{AB}}t_B+\overline{t_B}t_{A^{-1}B}+\overline{t_A}(1-|t_B|^2)$.

{\rm vii)}  $t_{ABA^2B^2}= t_{[A,B]} + t_{AB}t_{A^2B^2} - t_{AB}\overline{t_{AB}}$.

\end{proposition} 
\begin{proof}
The first identity follows from the definition of $\SU.$ The next two identities 
follow from the fact that the characteristic polynomial of $A$ has the 
form $A^3-t_AA^2+\overline{t_A}A-I$ (see the proof of Theorem 6.2.4 in \cite{goldman}). 

 Notice that by the Cayley-Hamilton theorem we have $A^3-t_AA^2+\overline{t_A}A-I=0.$ 
Now by multiplying this equality from the right by $A^{-1}B$ and then taking 
the trace, we obtain iv).

By multiplying the Cayley-Hamilton identity for $A$ by $A^{-1}B^2$ from the right and using previous identities we get v).

To prove vi) we will combine two equalities. The first one is obtained by multiplying the Cayley-Hamilton identity for $AB$ from the right by $(AB)^{-1}B^{-2}$:
$$ABAB^{-1}- t_{AB}AB^{-1}+\overline{t_{AB}}B^{-2}-B^{-1}A^{-1}B^{-2}=0.$$

The second one is obtained by multiplying the Cayley-Hamilton identity for $B$ from the left by $(AB)^{-1}B^{-2}$ :
$$ B^{-1}A^{-1}B-t_BB^{-1}A^{-1}+\overline{t_{B}}B^{-1}A^{-1}B^{-1}- B^{-1}A^{-1}B^{-2}=0.$$

It is not hard to see that when combining these two equalities and simplifying things by the previously proved identities we get the result.

The last identity can be obtained by multiplying the Cayley-Hamilton identity for $AB$ by
$(AB)^{-1}B^{-1}AB^2.$
\end{proof}  

We now state some algebraic results on the algebra of invariants of matrices. Let $\C[M_n^{\oplus m}]$ be the coordinate ring for the space of $m$-tuples 
of $n\times n$ matrices $(A_k=(a_{ij}^k))_{k=1,...,m},$ i.e., $\C[M_n^{\oplus m}]:=\C[a_{ij}^k| 1\le i,j\le n,1\le k \le m].$
Consider the action of $GL_n:={\rm GL}(n,\C)$ by simultaneous conjugation of $m$ matrices. Algebraists are interested in the algebra of invariants 
$C_{n,m}:=\C[M_n^{\oplus m}]^{GL_n}$

The following result of  
Teranishi \cite{teranishi} will be useful for us: 
The algebra $C_{3,2}$ of invariants of two  matrices $X, Y$ in ${\rm GL}(3,\C),$ is generated by :  
$$t_{X}, \, t_{Y },\,  t_{X^2}, \, t_{XY }, \, t_{Y^2},\, 
t_{X^3}, \, t_{X^2Y }, \, t_{XY^2}, \, t_{Y^3}, \, t_{X^2Y^2}, \, t_{X^2Y^2XY}$$  

This result means that the trace of any word in $X,Y$ can be expressed as a polynomial in the eleven traces above. When working with the group $\SU$
 we can reduce the number of generators greatly by  using Proposition \ref{trid}. We get the following :
\begin{proposition} \label{5}Let $A$ and $B$ be a pair of matrices in $\SU$ then the trace of any word in $A, B$ can be written as a polynomial of the following 
variables and their complex conjugates:
$$t_A, \, t_B, \, t_{AB}, \, t_{A^{-1}B}, \, t_{[A,B]}.$$
\end{proposition}  

Since the real dimension of $\SU$ is $8$, the real dimension of the representation space of the free group of rank 2 should also be $8.$ Therefore there should be a relation among these
 $5$ traces. Fortunately, this relation has been computed in \cite{ads,nakamoto} as the defining relation for the algebra of invariants. In particular, it has been 
shown that the algebra of invariants of two matrices in ${\rm GL}(3,\C)$ is defined  by a single relation which expresses  
$t_{X^2Y^2XY}$ as a solution of a quadratic equation whose coefficients are polynomials in the other ten traces. After plugging our variables into the 
formula in Theorem 1.2 of \cite{ads} and simplifying by MAPLE, we get the following result.
\begin{proposition} Let $A$ and $B$ be two matrices in $\SU.$ If we denote $t_A, t_B, t_{AB}, t_{A^{-1}B}$ by $a,b,c,d$ respectively,  then the following identities hold: 
\begin{eqnarray*}
\Re (t_{[A,B]})&=& \frac{1}{2}(|ab|^2
 + |a|^2 + |b|^2 + |c|^2 +|d|^2-ab\overline { c
 } -\overline {ab}c-a\overline { b } d-
\overline { a } b\overline { d } -3). \\
\Im(t_{[A,B]})^2 &=&-\frac{1}{4}(|ab|^2
 - |a|^2 - |b|^2 + |c|^2 +|d|^2-ab\overline { c
 } -\overline {ab}c-a\overline { b } d-
\overline { a } b\overline { d } )^2   \\
&& +2\Re( -a^3|b|^2+a^2\overline{b}^2\overline{d}+a^2b^2c-a|b|^2\overline{d}c-|a|^2b^3-|a|^2bcd+a^2\overline{c}d  \\
&& +a^2\overline{b}c+a^2\overline{d}b+ab^2d-2abc^2+acd^2
+\overline{b}\overline{d}c^2+b^2c\overline{a}-2bd^2\overline{a}+c^2d\overline{a} \\
&& +a^3+\frac{3}{2}ab\overline{c}+\frac{3}{2}a\overline{b}d-3ac\overline{d}+b^3+b^2\overline{cd}-3bcd+c^3+d^3+d^2\overline{bc} ) \\
&&+\frac{5}{2}|ab|^2+|cd|^2-\frac{9}{2}(|a|^2 + |b|^2 + |c|^2 +|d|^2)+\frac{27}{4}.  
\end{eqnarray*}

\end{proposition}      
\section{ Parameterization of the representation space}    
In this section we will show that the traces of certain elements give a coordinate system for the irreducible part of the representation space
 of the  Brieskorn homology sphere. Furthermore we also show how to determine the constraint region for the coordinates. 

\begin{theorem} Two irreducible representations $\rho, \rho':\fg(\bhs)\longrightarrow \SU$ are conjugate if and only if
the image under  $\rho$ and $\rho'$ of 
each $x,y,h$ are conjugate and  
satisfy the relations
$$t_{\rho(xy)}=t_{\rho'(xy)},
\quad
t_{\rho(x^{-1}y)}=t_{\rho'(x^{-1}y)}, 
\quad
\Im(t_{\rho([x,y])})=\Im(t_{\rho'([x,y])}).$$  
\end{theorem}

\begin{proof} 
If $\rho$ and $\rho'$ are conjugate, then
the required relations are obviously satisfied.
On the contrary suppose the relations hold.
Since $\rho$ and $\rho'$ are irreducible, $\rho(h)$ and $\rho'(h)$ should be in the center $Z(\SU)$ of $\SU.$ Notice that the images of $x,y,z$ 
 under a representation are elliptic elements, 
and they are diagonalizable. Moreover, it also follows from the irreducibility
that either $\rho(x)$ or $\rho(y)$ has three distinct 
eigenvalues since otherwise $\rho$ 
would have a non-trivial invariant subspace by 
dimensional reason.
The same holds  for $\rho'$.
So, after conjugation, we may assume that $\rho(x)=\rho'(x)= \diag(e^{i\theta_1}, e^{i\theta_2}, e^{i\theta_3}),$ where 
$e^{i\theta_1}, e^{i\theta_2}, e^{i\theta_3}$ are distinct numbers and $\diag(a,b,c)$ denotes the diagonal matrix whose diagonal elements are $a, b, c.$ 

To prove the theorem, it is enough to show that we can conjugate $\rho(y)$ to $\rho'(y)$ by a diagonal matrix. 
To show this we prepare a small lemma.
\begin{lemma} Let $A=\diag(e^{i\theta_1}, e^{i\theta_2}, e^{i\theta_3})$, where 
$e^{i\theta_1}, e^{i\theta_2}, e^{i\theta_3}$ are three distinct numbers. Suppose that $B=(b_{ij})$ and $B'=(b'_{ij})$ are two $3\times 3$ matrices 
satisfying
$t_B=t_{B'}, t_{AB}=t_{AB'}, t_{A^{-1}B}=t_{A^{-1}B'}$.
Then the diagonal elements of $B$ and $B'$ are equal.
\end{lemma}
\begin{proof}
It follows from 
our assumption that the following equations hold:
$$
\left\{\begin{array}{llll}
(b_{11}-b'_{11})   & + \ (b_{22}-b'_{22})  & +\ (b_{33}-b'_{33})& =0\\
(b_{11}-b'_{11})e^{i\theta_1}  & +\ (b_{22}-b'_{22})e^{i\theta_2}  & +\ (b_{33}-b'_{33})e^{i\theta_3}& =0\\
(b_{11}-b'_{11})e^{-i\theta_1}  & +\ (b_{22}-b'_{22})e^{-i\theta_2}  & +\ (b_{33}-b'_{33})e^{-i\theta_3}& =0
\end{array}\right.
$$ 
Consider this as a system of linear equations in $(b_{ii}-b'_{ii}).$ 
Since the determinant  
$$
\det\left(\begin{array}{ccc}
1&1&1\\
e^{i\theta_1}&e^{i\theta_2}&e^{i\theta_3}\\
e^{-i\theta_1}&e^{-i\theta_2}&e^{-i\theta_3}
\end{array}\right)-
\frac{(1-e^{i(\theta_1-\theta_2)})(1-e^{i(\theta_2-\theta_3)})(1-e^{i(\theta_3-\theta_1)})}
{e^{i(\theta_1+\theta_2+\theta_3)}}
$$
is not zero,
  we get the conclusion of the lemma.
\end{proof}      
Now come back to the proof of our theorem, suppose that $\rho(y)=(y_{ij})$ and $\rho'(y)=(y'_{ij}).$ 
From our assumption and Proposition 2.3, 
we have 
$t_{\rho(w)}=t_{\rho'(w)}$  for every word $w(x,y).$ Applying Lemma 3.2 for $A=\rho(x),$ $B=\rho(y)$ and $B'=\rho'(y),$ 
we obtain        
$$ (*)  \qquad \qquad \qquad \qquad y_{i,i}=y'_{i,i} \qquad (i=1,2,3).$$  
Applying Lemma 3.2 again 
for
$A=\rho(x),$ $B=\rho([x,y])$ and $B'=\rho'([x,y]),$ we conclude that the corresponding diagonal elements of $\rho([x,y])$
 and $\rho'([x,y])$ are equal. For the first diagonal element, we have:
$$  |y_{11}|^2 + |y_{12}|^2e^{i(\theta_1-\theta_2)} -|y_{13}|^2e^{i(\theta_1-\theta_3)} = |y'_{11}|^2 + |y'_{12}|^2e^{i(\theta_1-\theta_2)} -
|y'_{13}|^2e^{i(\theta_1-\theta_3)}.$$ 
Using the fact that $\rho(x), \rho(y)$ belong to $\SU$ and $y_{11}=y'_{11},$ we get the following 
equations:
$$
\left\{\begin{array}{lll} 
(|y_{12}|^2 - |y'_{12}|^2) &- (|y_{13}|^2 - |y'_{13}|^2)& =0\\
(|y_{12}|^2 - |y'_{1,2}|^2)e^{i(\theta_1-\theta_2)} &- (|y_{1,3}|^2 - |y'_{1,3}|^2)e^{i(\theta_1-\theta_3)}&=0.
\end{array}\right.
$$
From these equations, it follows that $|y_{12}|= |y'_{12}|$ and $|y_{13}|=|y'_{13}|$. Arguing similarly for other diagonal elements of $\rho([x,y])$
 and $\rho'([x,y]),$ we obtain that  
$$ (**)  \qquad \qquad\qquad \qquad |y_{ij}|= |y'_{ij}| \qquad ( i\ne j).$$ 

Applying Lemma 3.2 one more time 
for
 $A=\rho(x),$ $B=\rho(y^2)$ and $B'=\rho'(y^2),$  it follows that the corresponding diagonal elements of $\rho(y^2)$
 and $\rho'(y^2)$ are equal. Combining with $(*)$, we 
obtain the following equalities :
$$ (***) \qquad \qquad  \qquad \qquad y_{ij}y_{ji}=y'_{ij}y'_{ji} \qquad ( i\ne j). $$
Now consider three pairs $(y_{ij},y_{ji})$ for $i < j.$ By the irreducibility of $\rho,$ at least two pairs are not equal to $(0,0).$ Without loss of 
generality, we may assume that, say, $y_{12}\ne 0$ and $y_{31}\ne 0.$ By conjugate $\rho$ by $\diag(e^{i\phi_1}, e^{i\phi_2}, e^{i\phi_3})$ for appropriate
 values of $\phi_i$ and using $(**),$ we may assume that $y_{12}=y'_{12}$ and $y_{31}=y'_{31}.$ Furthermore, using $(***),$ we get that 
$y_{21}=y'_{21}$ and $y_{13}=y'_{13}.$ Moreover, since $\rho(y)$ and $\rho'(y)$ are in $\SU,$ we 
obtain
$$y_{11}\overline{y_{12}}+y_{21}\overline{y_{22}}-y_{31}\overline{y_{32}}=0,\qquad
y'_{11}\overline{y'_{12}}+y'_{21}\overline{y'_{22}}-y'_{31}\overline{y'_{32}}
=0.$$
It follows that $y_{32}=y'_{32}.$ By a similar argument, we also get $y_{23}=y'_{23}$ and thus our theorem is proved.      
\end{proof}
To describe the representation space, for each $h\in Z(\SU), x=\diag(\lambda_1, \lambda_2, \lambda_3), y= P \diag(\mu_1, \mu_2, \mu_3) P^{-1}, P\in \SU$ such that
$x^ph^a=y^qh^b=I,$ we need to answer 
the following
two questions:
 
- Does there exist $P$ such that $z=(xy)^{-1}$ satisfies $z^rh^c=I$?

- What are the possible values of $t_{x^{-1}y}$? 

In other words, we need to find the image of the following map in terms of $\lambda=(\lambda_1, \lambda_2, \lambda_3)$ and $\mu= (\mu_1, \mu_2, \mu_3)$ :
$$ \begin{array}{lll}   
\Phi_{\lambda,\mu} &:& \SU \longrightarrow  \C^2 \\
&& P \qquad  \mapsto  (t_{xy},t_{x^{-1}y}),
\end{array}$$
where $x=\diag(\lambda_1, \lambda_2, \lambda_3)$ and $y= P \diag(\mu_1, \mu_2, \mu_3) P^{-1}.$

If we write 
 $P=(p_{ij}),$, then we have
 $$
P^{-1}= \left(\begin{array}{ccc}
\overline{p_{11}}&\overline{p_{21}}&-\overline{p_{31}}\\
\overline{p_{12}}&\overline{p_{22}}&-\overline{p_{32}}\\
-\overline{p_{13}}&-\overline{p_{23}}&\overline{p_{33}}
\end{array}\right).$$
Let us denote by $\hat P$ for the matrix
$$\hat P=\left(\begin{array}{ccc}
\ \ |{p_{11}}|^2&\ |{p_{12}}|^2&-|{p_{13}}|^2\\
\ \ |{p_{21}}|^2&\ |{p_{22}}|^2&-|{p_{23}}|^2\\
-|{p_{31}}|^2&-|{p_{32}}|^2&\ \ |{p_{33}}|^2
\end{array}\right).$$
Then we have
$$t_{xy}=\left(\lambda_1, \lambda_2, \lambda_3\right)\hat P \left(\mu_1, \mu_2, \mu_3\right)^T, 
\qquad
t_{x^{-1}y}= 
\left(\frac{1}{\lambda_1},\frac{1}{\lambda_2},\frac{1}{ \lambda_3} \right)\hat P \left(\mu_1, \mu_2, \mu_3\right)^T.$$ 

Let $\mathcal D$ to be the set of $3\times 3$ matrices $M$ 
such that there exists $P=(p_{ij})\in \SU$ satisfying $M=\hat P.$
 An explicit description of $\mathcal D$ in the following lemma will help us to find the image of $\Phi_{\lambda,\mu}$ in practice.
\begin{lemma} 
Let $M$ be the matrix
$$M=\left(\begin{array}{ccc}
\ \ {m_{11}}&\ {m_{12}} &-{m_{13}}\\
\ \ {m_{21}}&\ {m_{22}}&-{m_{23}}\\
-{m_{31}}&-{m_{32}}&\ \ {m_{33}}
\end{array}\right)
$$ 
such that $m_{ij}\ge 0$ and the sum of every row or column is $1.$
Then $M$ is an element of $\mathcal D$ if and only if 
the following triangle inequalities holds: 
$$\sqrt{m_{1k}m_{2k}}\le \sum_{i\ne k} \sqrt{m_{1i}m_{2i}}
\qquad (k=1,2,3).
$$
\end{lemma}
\begin{proof} 
We first show the if part: 
If $M\in \mathcal D$ then there exist $\theta_{ij}$ such that the matrix $(\sqrt{m_{ij}}e^{i\theta_{ij}})$ belongs to $\SU.$,
and hence we have
$$\sqrt{m_{11}m_{21}}e^{i(\theta_{11}-\theta_{21})} + \sqrt{m_{12}m_{22}}e^{i(\theta_{12}-\theta_{22})} - \sqrt{m_{13}m_{23}}e^{i(\theta_{13}-\theta_{23})}=0.$$
It implies that 
the three numbers $\sqrt{m_{1i}m_{2i}}$ $(i=1,2,3)$ must satisfy the triangle inequalities and 
the lemma follows.

We next show the ``only if" part:
Now suppose that  three numbers $\sqrt{m_{1i}m_{2i}}$  
$(i=1,2,3)$ satisfy the triangle inequalities. 
Then
 there exist angles $\theta_{ij}$ 
satisfying
$$\sqrt{m_{11}m_{21}}e^{i(\theta_{11}-\theta_{21})} + \sqrt{m_{12}m_{22}}e^{i(\theta_{12}-\theta_{22})} - \sqrt{m_{13}m_{23}}e^{i(\theta_{13}-\theta_{23})}=0.$$
Put $p_{ij}=\sqrt{m_{ij}}e^{i\theta_{ij}}$ $(i=1,2, j=1,2,3)$, 
then
we get the first two rows of the matrix $P.$ Let $v=(p_{31},p_{32},p_{33})$ be the vector 
which is orthogonal to these two rows
with respect to
the indefinite inner product 
$\langle , \rangle_{2,1}$,  
and also satisfies $\langle v, v\rangle_{2,1}=-1.$ Then it is not hard to check 
 $P=(p_{ij})\in \SU.$ and $M=\hat P.$ 
\end{proof}

\section{examples}
{\bf Example 1.} The first example is the manifold $\Sigma(2,3,11).$ Its fundamental group 
has the following presentation.
$$\fg(\Sigma(2,3,11))= \langle
x,y,z,h| \ h\ \text{central},\ x^2h^{-1}=y^3h=z^{11}h^2=xyz=1
\rangle.
$$
In this example the irreducible representation space consists of isolated points. 

For each $h=\diag(\epsilon,\epsilon,\epsilon),$ 
with $\epsilon^3=1,$,
 we look for $\lambda$ and $\mu$ 
satisfying
$$\lambda_i^2=\epsilon, \ \  \mu_i^3=\epsilon^{-1} \quad (i=1,2,3), 
\qquad
\lambda_1\lambda_2\lambda_3=1,
\qquad
\mu_1\mu_2\mu_3=1$$ 
 such that the image of $\Phi_{\lambda,\mu}$ 
contains a point whose first coordinate 
is of 
 the form $e^{i\theta_1}+ e^{i\theta_2}+ e^{i\theta_3}$
 with $\theta_1+\theta_2+\theta_3=2k\pi$
and
$ e^{11i\theta_i}=\epsilon^2$  for all i.
 
 A small computer search tells us that there are five irreducible representations into $\SU$, all corresponding to the case 
$\rho(h)=I.$ The parameters of these representations 
 are given below. Here we use $\sim$ to denote the conjugacy relation.
 
 1) $\rho(x) \sim \diag(1,-1,-1), \quad \rho(y) \sim 
\diag(1,e^{4\pi i/3},e^{2\pi i/3}),$

 $t_{\rho(xy)}=t_{\rho(x^{-1}y)}= e^{10\pi i/11}+e^{16\pi i/11}+e^{18\pi i/11},\ \Im(t_{\rho([x,y])})=0.$
 
 2) $\rho(x) \sim \diag(1,-1,-1), \quad \rho(y) \sim 
\diag(e^{2\pi i/3},1,e^{4\pi i/3}),$

 $t_{\rho(xy)}=t_{\rho(x^{-1}y)}= e^{4\pi i/11}+e^{6\pi i/11}+e^{12\pi i/11},\ \Im(t_{\rho([x,y])})=0.$

3) $\rho(x) \sim \diag(-1,-1,1), \quad \rho(y) \sim 
\diag(1,e^{4\pi i/3},e^{2\pi i/3}),$

 $t_{\rho(xy)}=t_{\rho(x^{-1}y)}= e^{4\pi i/11}+e^{8\pi i/11}+e^{10\pi i/11},\ \Im(t_{\rho([x,y])})=0.$

4) $\rho(x) \sim \diag(-1,-1,1), \quad \rho(y) \sim 
\diag(e^{2\pi i/3},1,e^{4\pi i/3}),$

 $t_{\rho(xy)}=t_{\rho(x^{-1}y)}= e^{12\pi i/11}+e^{14\pi i/11}+e^{18\pi i/11},\ \Im(t_{\rho([x,y])})=0.$ 

5) $\rho(x) \sim \diag(-1,-1,1), \quad \rho(y) \sim 
\diag(e^{2\pi i/3},e^{4\pi i/3},1),$

 $t_{\rho(xy)}=t_{\rho(x^{-1}y)}= 1+2\cos(2\pi/11),\ \Im(t_{\rho([x,y])})=0.$
 
 It is no surprise that   $t_{\rho(xy)}=t_{\rho(x^{-1}y)}$ and
 $\Im(t_{\rho([x,y])})=0$ in all the cases since 
$\rho(x)^2=I.$ We can easily check that
these representations give 5 distinct irreducible representations when considered as elements of $\RR^*_{\PU}(\Sigma(2,3,11)).$

The Toledo invariant for representations of the fundamental group of an oriented surface into ${{\rm PU}(p,1)}$ is defined in \cite{toledo}.
In \cite{krebs, krebs2}, 
M.~Krebs defines the Toledo invariant for orbifold fundamental groups and uses it to obtain a lower bound for the number of connected components of    
the $\PU$ representation space. In particular, it is shown in \cite{krebs} that $\RR^*_{\PU}(\Sigma(2,3,11))$ has at least 5 connected components. So in this case the bound obtained by using the Toledo invariant is sharp.

{\bf Example 2.}  Our next example is the manifold $\Sigma(2,3,13)$ which has the fundamental group :
$$\fg(\Sigma(2,3,13))= \langle x,y,z,h| \ h\ \text{central},\ x^2h=y^3h^{-1}=z^{13}h^{-2}=xyz=1\rangle.$$
 A similar computer search as in the previous example shows that the irreducible representation space
 $\RR^*_{\SU}(\Sigma(2,3,13))$ consists of $8$ isolated points. In the following, we list the parameters of these representations. Note that 
$\rho(h)=I$ in all the cases.
 
 1) $\rho(x) \sim \diag(-1,1,-1),  \quad \rho(y) \sim 
\diag(e^{2\pi i/3},e^{4\pi i/3},1),$

 $t_{\rho(xy)}=t_{\rho(x^{-1}y)}= e^{4\pi i/13}+e^{10\pi i/13}+e^{12\pi i/13},\ \Im(t_{\rho([x,y])})=0.$
 
 2) $\rho(x) \sim \diag(-1,1,-1), \quad \rho(y) \sim 
\diag(e^{2\pi i/3},e^{4\pi i/3},1),$

 $t_{\rho(xy)}=t_{\rho(x^{-1}y)}= e^{14\pi i/13}+e^{16\pi i/13}+e^{22\pi i/13},\ \Im(t_{\rho([x,y])})=0.$

3) $\rho(x) \sim \diag(-1,1,-1), \quad \rho(y) \sim 
\diag(1,e^{4\pi i/3},e^{2\pi i/3}),$

 $t_{\rho(xy)}=t_{\rho(x^{-1}y)}= e^{6\pi i/13}+e^{22\pi i/13}+e^{24\pi i/13},\ \Im(t_{\rho([x,y])})=0.$

4) $\rho(x) \sim \diag(-1,1,-1), \quad \rho(y) \sim 
\diag(1,e^{2\pi i/3},e^{4\pi i/3}),$

 $t_{\rho(xy)}=t_{\rho(x^{-1}y)}= e^{2\pi i/13}+e^{4\pi i/13}+e^{20\pi i/13},\ \Im(t_{\rho([x,y])})=0.$

5) $\rho(x) \sim \diag(-1,-1,1), \quad \rho(y) \sim 
\diag(1,e^{4\pi i/3},e^{2\pi i/3}),$

 $t_{\rho(xy)}=t_{\rho(x^{-1}y)}= e^{6\pi i/13}+e^{8\pi i/13}+e^{12\pi i/13},\ \Im(t_{\rho([x,y])})=0.$
 
6)$\rho(x) \sim \diag(-1,-1,1), \quad \rho(y) \sim 
\diag(e^{2\pi i/3},e^{4\pi i/3},1),$

 $t_{\rho(xy)}=t_{\rho(x^{-1}y)}= 1+2\cos(2\pi/13),\ \Im(t_{\rho([x,y])})=0.$

7)$\rho(x) \sim \diag(-1,-1,1), \quad \rho(y) \sim 
\diag(e^{2\pi i/3},e^{4\pi i/3},1),$

 $t_{\rho(xy)}=t_{\rho(x^{-1}y)}= 1+2\cos(4\pi/13),\ \Im(t_{\rho([x,y])})=0.$
 
8) $\rho(x) \sim \diag(-1,-1,1), \quad \rho(y) \sim 
\diag(1,e^{2\pi i/3},e^{4\pi i/3}),$

 $t_{\rho(xy)}=t_{\rho(x^{-1}y)}= e^{14\pi i/13}+e^{18\pi i/13}+e^{20\pi i/13},\ \Im(t_{\rho([x,y])})=0.$ 
 
 These representations give us $8$ distinct points of $\RR^*_{\PU}(\Sigma(2,3,13)).$ According to \cite{krebs2}(Appendix), in this case there are $7$ distinct values of the orbifold Toledo invariant. So the    
orbifold Toledo invariant does not distinguish the connected components of $\RR^*_{\PU}(\Sigma(2,3,13)).$

{\bf Acknowledgment.}  
The author would like to thank Professor Takashi Tsuboi for advice and hospitality during the time in Tokyo. The author is grateful to his former thesis advisor, Daniel
Ruberman, for informing him of the work of 
M.~Krebs and for continuous support. We express our sincere thanks to the anonymous referee for pointing out several inaccuracies in the earlier version of this paper.  
 

\end{document}